\documentclass{amsart}

\newtheorem{theorem}{Theorem}
\newcommand{\bt}{\begin{theorem}}
\newcommand{\et}{\end{theorem}}
\newtheorem{lemma}{Lemma}
\newcommand{\bl}{\begin{lemma}}
\newcommand{\el}{\end{lemma}}
\newtheorem{corollary}{Corollary}
\newcommand{\bc}{\begin{corollary}}
\newcommand{\ec}{\end{corollary}}
\newtheorem{problem}{Problem}
\newcommand{\bprob}{\begin{problem}}
\newcommand{\eprob}{\end{problem}}

\newtheorem{example}{Example}
\newcommand{\bex}{\begin{example}}
\newcommand{\eex}{\end{example}}

\newcommand{\beq}{\begin{equation}}
\newcommand{\eeq}{\end{equation}}
\newcommand{\benum}{\begin{enumerate}}
\newcommand{\eenum}{\end{enumerate}}
\newcommand{\N}{\ensuremath{ \mathbf N }}
\newcommand{\Z}{\ensuremath{\mathbf Z}}

\newcommand{\mcf}{\ensuremath{ \mathcal F}}

\newcommand{\bsmallmat}{\left(\begin{smallmatrix}}
\newcommand{\esmallmat}{\end{smallmatrix}\right)}

\DeclareMathOperator{\id}{id}
\DeclareMathOperator{\orbit}{\text{orbit}}
\newcommand{\bmat}{\left(\begin{matrix}}
\newcommand{\emat}{\end{matrix}\right)}
\DeclareMathOperator{\qand}{\quad\text{and}\quad}
\DeclareMathOperator{\qqand}{\qquad\text{and}\qquad}

\DeclareMathOperator{\Perm}{\text{Perm}}

\DeclareMathOperator{\lcm}{\text{lcm}}

\usepackage[all]{xy}
\usepackage{amssymb,latexsym}

\title[Fixed points of powers of permutations]
{Arithmetic functions and fixed points of powers of permutations}
\author{Melvyn B.  Nathanson}
\address{Department of Mathematics\\Lehman College (CUNY)\\Bronx, NY 10468} 
\email{melvyn.nathanson@lehman.cuny.edu} 
\date{\today}

\subjclass[2010]{11N56, 20B05, 20B07, 20B10, 20F69}

\keywords{Permutations, fixed points, cycle structure, arithmetic functions, group theory.}

\begin{document}

\begin{abstract}
Let $\sigma$ be a permutation of a nonempty finite or countably infinite set $X$  
and let $F_X\left( \sigma^k\right)$ count the number of fixed points of 
the $k$th power of $\sigma$.
This paper explains how the arithmetic function 
$k \mapsto \left(F_X\left( \sigma^k\right) \right)_{k=1}^{\infty}$ 
determines the conjugacy class of the permutation $\sigma$, 
constructs an algorithm to compute the conjugacy class from the fixed point counting function 
$F_X\left( \sigma^k\right)$, and describes the arithmetic functions that are fixed point 
counting functions of permutations. 
\end{abstract}

\maketitle

\section{Fixed point counting functions} 
Let $X$ be a nonempty finite or countably infinite set 
and let $|X|$ denote the cardinality of $X$.
Let $\Perm(X)$ be the group of permutations of $X$. 
The identity permutation is $\id_X$.  
We write $S_n = \Perm(X)$ if $X = \{1,2,3,\ldots, n\}$.
Let \Z\ be the set of integers, $\N = \{ 1,2, 3,\ldots \}$ the set of positive integers, 
and $\N_0 = \N \cup \{ 0\} = \{0, 1, 2, 3, \ldots\}$ the set of nonnegative integers. 

Let $\ell \in \N$ and let  $\{x_1,\ldots, x_{\ell}\}$ be a finite set of distinct elements of $X$.   
The permutation $\gamma = (x_1,\ldots, x_{\ell}) \in \Perm(X)$
 defined by 
 \[
 \gamma(x_i) = x_{i+1} \text{ for } i \in \{1,\ldots, \ell -1\},
\qquad 
\gamma(x_{\ell}) = x_1,
\]
and 
\[
\gamma(x) = x \text{ for all } x \in X\setminus \{x_1,\ldots, x_{\ell} \}  
\]
is a \index{cycle}\index{cycle!lenght}\emph{cycle of length $\ell$}.  
Let $\{x_i:i\in \Z\}$ be an infinite set of distinct elements of $X$. 
The permutation $\gamma = (\ldots, x_{-2}, x_{-1}, x_0, x_1, x_2, \ldots, ) \in \Perm(X)$ 
defined by $\gamma(x_i) = x_{i+1}$ for all $i \in \Z$ is a cycle of infinite length, that is, 
\emph{cycle of  length $\infty$ }. 

Every permutation $\sigma \in \Perm(X)$ has a unique decomposition  
as a product of pairwise disjoint cycles of finite 
or infinite length.  
This decomposition is called the \emph{cyclic representation}\index{cyclic representation} 
of the permutation.  A permutation $\sigma$ \emph{contains the cycle $\gamma$} 
if $\gamma$ appears in the cyclic representation of $\sigma$.  
For example, the permutation $\sigma_0 \in \Perm(\N)$ defined by 
\[
\sigma_0 = (1) (2,3) (4,5,6) (7,8,9,10) (11, 12, 13, 14, 15)  \cdots 
\]
contains one cycle of length $\ell$ for all $\ell \in \N$. 
The cycle $(7,9)$ is contained in $\sigma_0^2$ but not in $\sigma_0$.

For all $\ell \in \N \cup \{\infty\}$, the 
\emph{cycle counting function}\index{cycle counting function} 
 $C_{\sigma}(\ell)$ counts the number of cycles of length $\ell$ contained in $\sigma$. 
Permutations $\sigma$ and $\tau$ in $\Perm(X)$ 
are  \emph{conjugate}\index{conjugate permutations}  if
$C_{\sigma}(\ell) = C_{\tau}(\ell)$ for all $\ell \in \N \cup \{\infty\}$.  
Equivalently, $\sigma$ and $\tau$ are conjugate if,
in their representations as products of pairwise disjoint cycles, 
the number of cycles of length $\ell$ in $\sigma$ equals the number of 
cycles of length $\ell$ in $\tau$  for all $\ell \in \N \cup \{\infty\}$.

The permutation $\sigma \in \Perm(X)$ is of  
\emph{finite type}\index{finite type}\index{permutation!finite type} 
if it is the product of finitely many or infinitely many pairwise disjoint cycles of finite length.  
Equivalently, 
$\sigma$ is of finite type if it contains no cycle of infinite length, 
that is, if $C_{\sigma}(\infty) = 0$.  
A permutation $\sigma \in \Perm(X)$ is of  
\emph{infinite type}\index{infinite type}\index{permutation!infinite type}
if it is the nonempty product of cycles of infinite length 
or, equivalently, if $C_{\sigma}(\infty) \geq 1$ and $C_{\sigma}(\ell) = 0$ for all $\ell \in \N$. 

The  permutation $\sigma$ has \emph{bounded finite cycle length} or, simply, is 
\emph{bounded}\index{bounded permutation}\index{permutation!bounded}
if there is a positive integer $\ell^*$ such that the length of every finite cycle 
in $\sigma$ is at most $\ell^*$, and so $C_{\sigma}(\ell) = 0$ for all $\ell \in \N$ 
with $\ell > \ell^*$.  
A permutation has 
\emph{unbounded finite cycle length}\index{unbounded permutation}\index{permutation!unbounded} 
or is \emph{unbounded} if it  contains finite cycles of arbitrarily large length.  
A permutation of bounded or unbounded  finite cycle length may 
also contain permutations of infinite length.  

The permutation $\sigma$ has \emph{finite multiplicity}\index{finite multiplicity} 
if it contains only finitely many cycles of length $\ell$ for all $\ell \in \N$, that is, 
if $C_{\sigma}(\ell) < \infty$ for all $\ell  \in \N$. 
The permutation $\sigma$ has \emph{infinite multiplicity}\index{infinite multiplicity} 
if $C_{\sigma}(\ell) = \infty$ for some $\ell  \in \N$. 
A nonempty set $X$ is finite if and only if every permutation in $\Perm(X)$ 
has finite type and  finite multiplicity and is bounded. 
A permutation has finite order in the group $\Perm(X)$  if and only if it is of finite type and bounded.

An \emph{arithmetic function}\index{arithmetic function} is a function 
whose domain is the set  $\N$ of positive integers. 
The \emph{zero arithmetic function}\index{arithmetic function!zero} 
is the function $f(k) = 0$ for all $k \in \N$. 
An arithmetic function $f$  is \emph{nonzero}\index{arithmetic function!nonzero} 
if $f(k) \neq 0$ for some $k \in \N$. 
An arithmetic function $f$  is \emph{nonnegative}\index{arithmetic function!nonnegative} 
if $f(k) \geq 0$ for all $k \in \N$.

Every nonzero arithmetic function $C:\N \rightarrow \N_0 \cup \{\infty\}$ is the cycle counting 
function of a unique (up to conjugacy) permutation of finite type 
on a nonempty finite or countably infinite set. 

A \emph{fixed point}\index{fixed point} of a permutation $\sigma \in \Perm(X)$ 
is an element $x \in X$ such that $\sigma(x) = x$.  Thus, $x$ is a fixed point if and only if $(x)$ 
is a cycle of length one in $\sigma$. 
Note that every element of $X$ is a fixed point of $\sigma$ if and only if $\sigma = \id_X$. 

For every integer $k$, let $\sigma^k$ be the $k$th power of $\sigma$.
An element $x \in X$ is a fixed point of $\sigma^k$ if $\sigma^k(x) = x$.  
The \emph{fixed point counting function}\index{fixed point counting function}   
$k \mapsto F_X\left( \sigma^k \right)$  is the arithmetic function that counts the number 
of fixed points of $\sigma^k$ for $k \in \N$.  
For example, if $\gamma$ is a cycle of length $\ell$, then
\[
F_X(\gamma^k) = \begin{cases}
\ell & \text{if $k \equiv 0 \pmod{ \ell}$} \\ 
0 & \text{if $k \not\equiv 0 \pmod{\ell}$.} 
\end{cases}
\]

If $x  \in X$ and $\sigma^{\ell}(x) = x$, then for all $d \in \N$ 
we have $\sigma^{d\ell}(x) = x$  and so 
the fixed point counting function of a permutation must satisfy the following condition: 
\[
F_X\left( \sigma^{d\ell} \right) \geq F_X\left( \sigma^{\ell} \right).
\]

The number $F_X\left( \sigma^k \right)$ can be finite or infinite.  
If $\sigma$ contains infinitely many cycles of length $k$, 
then $F_X\left( \sigma^k \right) = \infty$.  We have $F_X\left( \sigma^k \right) = \infty$ 
if and only if $\sigma$ contains infinitely many cycles of length ${\ell}$ for some divisor ${\ell}$ of $k$. 
Thus, $F_X(\sigma^k) < \infty$ for all $k \in \N$ if and only if $\sigma$ has finite multiplicity. 

If $ (\ldots, x_{-2}, x_{-1},x_0,x_1,x_2,\ldots)$ is an infinite cycle in $\sigma$, 
then $\sigma^k(x_i) = x_{i+k} \neq x_i$ for all $k \in \N$. 
It follows that if $\sigma$ is a permutation of infinite type, 
then $F_X(\sigma^k) = 0$ for all $k \in \N$ and so 
the fixed point counting function of a permutation is the zero function 
if and only if the permutation is a product of infinite cycles. 
Equivalently, $F_X(\sigma^k) \neq 0$ for some $k \in \N$ if and only if $\sigma$ contains a cycle 
of finite length.

For $i, \ell \in \N$, let 
\[
\gamma_{i,\ell} = (i,i+1,i+2,\ldots, i+\ell-1)
\]
be  the cycle of length $\ell$  in $\Perm(\N)$.   
For $a,d \in \N$ and $i \in \Z$, let 
\[
x_i = \begin{cases}
a+(2i-1)d & \text{if $i \geq 1$} \\
a-2id & \text{if $i \leq 0$}
\end{cases}
\]
and consider the infinite cycle 
\[
\Gamma_{a,d} = (x_i)_{i \in \Z} = (\ldots, a+4d, a+2d, a, a+d, a+3d, a+5d, \ldots). 
\]
For example, 
\[
\Gamma_{3,3} = (x_i)_{i \in \Z} = (\ldots, 15, 9, 3, 6, 12, 18, \ldots). 
\]
Let $t_j = \sum_{i=1}^j i  = j(j+1)/2$ be the $j$th triangular number. 

Consider the following permutations of \N: 
\begin{align*}
\sigma_0 & = \prod_{j=1}^{\infty}  \gamma_{t_{j-1}+1,j} 
= (1) (2,3) (4,5,6) (7,8,9,10) (11, 12, 13, 14, 15)  \cdots \\
\sigma_1 & = \prod_{j=0}^{\infty} \gamma_{2j+1,2}  =  (1,2)(3,4)(5,6)(7,8)(9,10) \cdots      \\
\sigma_2 & =  \gamma_{1,4} \prod_{j=2}^{\infty} \gamma_{2j+1,2}  =  (1,2,3,4)(5,6)(7,8)(9,10) \cdots  \\
\sigma_3 & =   \prod_{j=0}^{\infty} \gamma_{6j+1,2}  \prod_{j=0}^{\infty} \gamma_{6j+3,4}  \\
& =  (1,2)(3,4,5,6)(7,8)(9,10,11,12)(13,14) \cdots \\
\sigma_4 & = \Gamma_{3,3}  \prod_{i=0}^{\infty} \gamma_{3i+1,2} \\ 
& = (\ldots, 15, 9, 3, 6, 12, 18, \ldots) (1,2) (4,5) (7,8) (10,11)  \cdots. 
\end{align*}
Computing the cycle counting functions for these permutations, we obtain  
\begin{align*}
C_{\sigma_0}(\ell) &  = 1  \quad \text{ for all $\ell \in \N$} \\
C_{\sigma_0}(\infty) & = C_{\sigma_1}(\infty)  = C_{\sigma_2}(\infty) = C_{\sigma_3}(\infty)  = 0  \qand  C_{\sigma_4}(\infty)  = 1\\
C_{\sigma_1}(\ell) & = C_{\sigma_2}(\ell)  = C_{\sigma_3}(\ell)  = C_{\sigma_4}(\ell)  = 0  \quad
\text{ for all $\ell \in \N\setminus \{2,4\}$} \\
C_{\sigma_1}(2) & = C_{\sigma_2}(2)  = C_{\sigma_3}(2)  = C_{\sigma_4}(2)   = \infty  \\ 
C_{\sigma_1}(4) & = C_{\sigma_4}(4)  = 0, \quad 
 C_{\sigma_2}(4)   =  1, \quad 
 C_{\sigma_3}(4)   = \infty. 
\end{align*}

The permutations $\sigma_0$, $\sigma_1$, $\sigma_2$, and $\sigma_3$ 
have finite type; the permutation $\sigma_4$ does not have finite type.   

The permutation $\sigma_0$ is unbounded.   
The permutations $\sigma_1$, $\sigma_2$, $\sigma_3$, and $\sigma_4$ are bounded.

The permutation $\sigma_0$ has finite multiplicity.   
The permutations $\sigma_1$, $\sigma_2$, $\sigma_3$, and $\sigma_4$ have infinite multiplicity.

The  permutations $\sigma_1$, $\sigma_2$, $\sigma_3$, $\sigma_4$ 
are pairwise nonconjugate, 
but have identical fixed point counting functions:  
\[
F_{\N}(\sigma_1^k) = F_{\N}(\sigma_2^k) = F_{\N}(\sigma_3^k) = F_{\N}(\sigma_4^k) = 
\begin{cases}
0 & \text{if $k$ is odd}\\
\infty & \text{if $k$ is even.} 
\end{cases}
\]

There is a large literature on fixed points of permutations. 
Some recent papers are Cameron and Cohen~\cite{came92}, 
Diaconis, Fulman,  and Guralnick~\cite{diac08}, and Ford~\cite{ford21}.  
It is a folklore result\footnote{By \emph{folklore} I mean that the result is known 
to the experts but I do not  have a reference.} 
in group theory that every permutation of a finite set is  determined, 
up to conjugacy, by its fixed point counting 
function.  The examples above show that this is not true for permutations of an infinite set.

We shall prove that the conjugacy class of a permutation $\sigma$ 
of a nonempty finite or countably infinite 
set is uniquely determined 
by its fixed point counting function if and only if the permutation has finite type 
and finite multiplicity. Moreover, 
we give a simple formula that explicitly computes the conjugacy class 
(that is,  computes the number $C_{\sigma}(\ell)$ 
of cycles of length $\ell$) 
in terms of $F_X(\sigma^k)$. 
We also describe the permutations that have periodic fixed point counting functions 
and we classify the arithmetic functions that are fixed point counting functions 
of permutations of finite type.

\section{Fixed points and the M\" obius function} 

\bl                                                                 \label{FixedPoint:lemma:disjointX} 
Let $(X_i)_{i\in I}$ be a family of pairwise disjoint nonempty sets and let 
$\sigma_i \in \Perm(X_i)$ for all $i \in I$.  Let $X = \bigcup_{i\in I} X_i$.
The permutation $\sigma:X \rightarrow X$ defined by 
\[
\sigma(x) = \sigma_i(x) \qquad \text{for all $x \in X_i$}
\]
satisfies the addition formula   
\[
F_X(\sigma^k) =   \sum_{i\in I} F_{X_i}(\sigma_i^k) \qquad \text{ for all $k \in \N$.}
\]
\el

\begin{proof}
It suffices to observe that if $\sigma^k(x) = x$, then there is a unique $i \in I$ such that $x \in X_i$ 
and $\sigma_i^k(x) = x$.  Conversely, if $x \in X_i$ and $\sigma_i^k(x) = x$, 
then $\sigma^k(x) = x$ and $x \notin X_j$ for all $j \neq i$. 
This completes the proof. 
\end{proof}

Let $X$ be a nonempty set.  For all $\sigma \in \Perm(X)$, define the sets 
\beq                                                                \label{FixedPoint:X-0}
\hat{X}_0 = \hat{X}_0(\sigma) 
= \{x \in X: \text{$x$ is an element in a cycle  of $\sigma$ of finite length} \}
\eeq
and 
\beq                                                         \label{FixedPoint:X-infty}
\hat{X}_{\infty} = \hat{X}_{\infty}(\sigma) 
= \{x \in X: \text{$x$ is an element in a cycle of $\sigma$ of infinite length} \}.
\eeq

\bt                                               \label{FixedPoint:theorem:PartitionFinite} 
Let $X$ be a nonempty set and let $\sigma \in \Perm(X)$.  
The sets 
$\hat{X}_0$ and $\hat{X}_{\infty}$ defined by~\eqref{FixedPoint:X-0}
and~\eqref{FixedPoint:X-infty} partition $X$:
\[
X = \hat{X}_0  \cup \hat{X}_{\infty}  \qqand  \hat{X}_0  \cap \hat{X}_{\infty} = \emptyset. 
\] 
Let $\hat{\sigma}_0$ be the restriction of $\sigma$ to $\hat{X}_0$ 
and let $\hat{\sigma}_{\infty}$ be the restriction of $\sigma$ to $\hat{X}_{\infty}$.  
Then  $\hat{\sigma}_0 \in \Perm(\hat{X}_0)$ and 
$\hat{\sigma}_{\infty} \in \Perm(\hat{X}_{\infty})$.
For all positive integers $k$,  
\[
F_X(\sigma^k) = F_{\hat{X}_0}\left(\hat{\sigma}_0^k\right) 
\] 
and 
\[
F_{ \hat{X}_{\infty} }\left(\hat{\sigma}_{\infty}^k \right) = 0.  
\] 
\et

\begin{proof} 
For all $x \in X$, the elements $x$,  $\sigma(x)$, and $\sigma^{-1}(x)$ 
belong to the same cycle in the cyclic decomposition of $\sigma$, and so  
$\hat{\sigma}_0 \in \Perm(\hat{X}_0)$ and 
$\hat{\sigma}_{\infty} \in \Perm(\hat{X}_{\infty})$.

If $x \in X$ is an element in a cycle of infinite length, then $\sigma^k(x) \neq x$ 
for all $k \in \N$ and so $F_{ \hat{X}_{\infty} }(\hat{\sigma}_{\infty}^k) = 0$ 
for all $k \in \N$.
The addition formula of Lemma~\ref{FixedPoint:lemma:disjointX} 
applied to the partition $X = \hat{X}_0  \cup \hat{X}_{\infty}$ implies 
\[
F_X(\sigma^k) = F_{ \hat{X}_{0} }(\hat{\sigma}_{0}^k) 
+  F_{ \hat{X}_{\infty} }(\hat{\sigma}_{\infty}^k) 
= F_{ \hat{X}_{0} }(\hat{\sigma}_{0}^k)
\]
for all $k \in \N$. 
This completes the proof. 
\end{proof}

The permutation $\hat{\sigma}_0 \in \Perm(\hat{X}_0 ) $ 
is called the \emph{finite part}\index{finite part} of $\sigma$ 
and the permutation $\hat{\sigma}_{\infty} \in \Perm(\hat{X}_{\infty} ) $ 
is called the \emph{infinite part}\index{infinite part} of $\sigma$.     
Theorem~\ref{FixedPoint:theorem:PartitionFinite} shows that the 
fixed point counting function of a permutation only ``sees'' the finite part of the permutation. 
From the fixed point counting function of a permutation one can hope only 
to reconstruct the finite part of the permutation. 

A  permutation of finite type may have infinite multiplicity, 
that is, for some positive integer $\ell$, the permutation 
contains infinitely many cycles of length $\ell$.  
The permutations $\sigma_1$,  $\sigma_2$, $\sigma_3$,  $\sigma_4$ 
constructed above have the same fixed point counting function 
but have infinite multiplicity and are not conjugate.  
Thus, the fixed point counting function cannot separate  
permutations of finite type and infinite multiplicity.  
The unique reconstruction of a permutation from its fixed point counting function  
is possible only for permutations $\sigma$ of finite type  and finite multiplicity. 
For such permutations there is a simple formula that inverts the fixed point counting function.  

The M\" obius function $\mu(k)$ is the arithmetic function 
\[
\mu(k) = \begin{cases}
1 & \text{if $k = 1$}\\
(-1)^r & \text{if $k$ is the product of $r$ distinct primes} \\
0 & \text{if $k$ is divisible by $p^2$ for some prime $p$.}
\end{cases}
\]
Let $f$ and $g$ be arithmetic functions.  
The M\" obius inversion theorem from elementary number theory 
(Nathanson~\cite[Theorem 6.14]{nath00})  
states that 
\[
g(\ell) =  \sum_{\substack{ k= 1 \\ k|\ell}}^{\ell} f(k)
\]
if and only if  
\[
f(\ell) =   \sum_{\substack{ k= 1 \\ k|\ell}}^{\ell}\mu\left( \frac{\ell}{k} \right)  g(k). 
\]

\bt              \label{FixedPoint:theorem:InversionFormula}
Let  $\sigma \in \Perm(X)$ be a  permutation of finite type.
For every positive integer $\ell$, 
\[
F_X(\sigma^{\ell}) = \sum_{\substack{ k= 1 \\ k|\ell}}^{\ell} k C_{\sigma}(k).
\]
Let  $\sigma \in \Perm(X)$ be a  permutation of finite type and finite multiplicity.
For every positive integer $\ell$, 
\[
C_{\sigma}(\ell) = \frac{1}{\ell}  \sum_{\substack{ k= 1 \\ k|\ell}}^{\ell} 
 \mu\left( \frac{\ell}{k} \right) F_X(\sigma^k)
\]
where $\mu(k)$ is the M\" obius function.  
Thus, the fixed point counting function  of a  permutation of finite type 
and finite multiplicity determines the conjugacy class of the permutation.    
\et

\begin{proof}
For each $x \in X$ there is a unique cycle $\gamma$ 
in the cyclic decomposition of $\sigma$ such that $\gamma$ contains $x$.   
Let $k$ be the length of the cycle $\gamma$. 
We have $\sigma^{\ell}(x) = x$ if and only if $\gamma^{\ell}(x) = x$ 
if and only if $\gamma^{\ell} = \id_X$ if and only if $k$ divides $\ell$
and so 
\beq         \label{FixedPoint:Mob}
F_X(\sigma^{\ell}) = \sum_{\substack{ k= 1 \\ k|\ell}}^{\ell} k C_{\sigma}(k).
\eeq
This formula is valid even if $\sigma$ has infinite multiplicity, that is, if 
$C_{\sigma}(k)  = \infty$ for some $k \in \N$.

The permutation $\sigma$ has finite multiplicity if and only if 
$C_{\sigma}(k) < \infty$ for all $k \in \N$. 
In this case, applying the M\" obius inversion formula to~\eqref{FixedPoint:Mob}, we obtain  
\[
\ell C_{\sigma}(\ell) =   \sum_{\substack{ k= 1 \\ k|\ell}}^{\ell}  \mu\left( \frac{\ell}{k} \right) F_X(\sigma^k). 
\]
This completes the proof. 
\end{proof} 

For example, for the permutation 
\[
\sigma_0  = (1) (2,3) (4,5,6) (7,8,9,10) (11, 12, 13, 14, 15)  \cdots \in \Perm(\N)
\] 
we have $C_{\sigma_0}(\ell) = 1$ for all $\ell \in \N$ and 
\[
F_X(\sigma^{\ell}) = \sum_{\substack{ k= 1 \\ k|\ell}}^{\ell} k C_{\sigma}(k)
= \sum_{\substack{ k= 1 \\ k|\ell}}^{\ell} k = s(k)
\]
where $s(k)$ is the usual ``sum of the divisors'' function. It follows that 
\[
\ell =  \sum_{\substack{ k= 1 \\ k|\ell}}^{\ell} 
 \mu\left( \frac{\ell}{k} \right) s(k). 
\]

Every permutation of a finite set is of finite type and multiplicity, and so 
Theorem~\ref{FixedPoint:theorem:InversionFormula} gives a quantitative 
expression of the folklore result cited in Section 1.

We also obtain a necessary and sufficient condition for an arithmetic function 
to be the fixed point counting function of a permutation of finite type and finite multiplicity.

\bt
An arithmetic function $f:\N \rightarrow \N_0$ is the fixed point counting function 
of a permutation of finite type and finite multiplicity if and only if, for all $\ell \in \N$, 
\[
 \sum_{\substack{ k= 1 \\ k|\ell}}^{\ell} 
 \mu\left( \frac{\ell}{k} \right) f(k) \geq 0
\]
and 
\[
 \sum_{\substack{ k= 1 \\ k|\ell}}^{\ell} 
 \mu\left( \frac{\ell}{k} \right) f(k) \equiv 0 \pmod{\ell}. 
\]
\et

\begin{proof}
Necessity follows from Theorem~\ref{FixedPoint:theorem:InversionFormula}. 
To prove sufficiency, we observe that, for all $\ell \in \N$, 
the arithmetic function $g:\N \rightarrow \N_0$ defined by 
\[
g(\ell) = \frac{1}{\ell} \sum_{\substack{ k= 1 \\ k|\ell}}^{\ell} 
 \mu\left( \frac{\ell}{k} \right) f(k)
\]
is a nonnegative integer and so there is a permutation $\sigma$ 
on a set $X$ whose 
cycle counting function is $C_{\sigma}(\ell)= g(\ell)$.  
By M\" obius inversion and Theorem~\ref{FixedPoint:theorem:InversionFormula},
\[
f(\ell) =  \sum_{\substack{ k= 1 \\ k|\ell}}^{\ell} k g(k) 
= \sum_{\substack{ k= 1 \\ k|\ell}}^{\ell} k C_{\sigma}(k) 
= F_X(\sigma^{\ell}). 
\]
This completes the proof.  
\end{proof}

\section{An alternative algorithm}

There is an alternative, algorithmic method to determine if 
an arithmetic function $f$ 
is the fixed point counting function of a permutation $\sigma$ 
of finite type and finite multiplicity
and also to compute the permutation $\sigma$, that is, to compute the 
cycle counting function $C_{\sigma}(k)$. 
The construction is based on three observations about 
 the fixed point counting function $F_X(\sigma^k)$ of a permutation $\sigma \in \Perm(X)$   
of finite type and finite multiplicity.
\benum
\item
The arithmetic function $k \mapsto F_X(\sigma^k)$ is integer-valued and nonnegative. 
\item 
If $x$ is a fixed point of $\sigma^{\ell}$ for some $\ell \in \N$, then 
$x$ is a fixed point of $\sigma^{d\ell}$ for all $d \in \N$ and so 
\[
F_X(\sigma^{d\ell}) - F_X(\sigma^{\ell})  \geq 0 \quad\text{for all $d \in \N$.}
\] 

\item
If $\ell$ is the least positive integer such that $ F_X(\sigma^{\ell} ) > 0$, 
then there exists $x \in X$ with $\sigma^{\ell}(x) = x$ 
and $x$ is an element in a cycle of length $\ell$.  
If $\sigma$ contains exactly $q$ cycles of length $\ell$, that is, if $C_{\sigma}(\ell) = q$, then 
\[
F_X(\sigma^{\ell}) = q\ell \equiv 0 \pmod{\ell}. 
\]
Let $W$ be the subset of $X$ consisting of the $q\ell$  elements of $X$ that 
appear in the $q$ cycles of length $\ell$, and let $Y = X\setminus W$.
The restriction of $\sigma$ to $Y$ is a permutation $\tau \in \Perm(Y)$ 
such that 
\[
F_Y(\tau^k) = 0 \qquad\text{ for all } k \leq \ell 
\] 
\[
F_Y(\tau^{d\ell}) = F_X(\sigma^{d\ell}) - q\ell
\qquad \text{ for all } d \in \N 
\]
and 
\[
F_Y(\tau^k) = F_X(\sigma^k) 
\qquad \text{ for all } k \not\equiv 0 \pmod{\ell}.
\]
\eenum

Let $f = f_1$ be a nonzero integer-valued arithmetic function.  
 We construct inductively a sequence of arithmetic functions 
$(f_j)_{j=1}^{\infty}$.   
\benum
\item
If $f_j(k) < 0$ for some $k \in \N$, then 
$f_{j+1}$ is the constant function -1. 
It follows that $f_i$ is the constant function -1 for all $i \geq j+1$. 

\item
If $f_j$ is the zero function, then $f_{j+1}$ is the zero function. 
It follows that $f_i$ is the zero function for all $i \geq j$. 

\item
If $f_j$ is a nonzero nonnegative arithmetic function, 
then $f_j(k) > 0$ for some $k \in \N$. Let $\ell_j$ be the smallest positive 
integer such that $f_j(\ell_j) > 0$. 
\benum
\item
If $f_j(\ell_j) \not\equiv 0 \pmod{\ell_j}$, 
then $f_{j+1}$  is the constant function -1.
\item
If $f_j(\ell_j) \equiv 0 \pmod{\ell_j}$, 
then $f(\ell_j) = q_j\ell_j$ for some positive integer $q_j$.
In this case, we construct the integer-valued arithmetic function $f_{j+1}$ as follows: 
\benum 
\item
If 
\[
f_{j}(d\ell_j ) - f_{j}(\ell_j )  = f_{j}(d\ell_j ) - q_j\ell_j < 0
\]
 for some $d \in \N$, 
then $f_{j+1}$ is the constant function -1.  
\item 
If 
\[
f_{j}(d\ell_j ) - f_{j}(\ell_j )  = f_{j}(d\ell_j ) - q_j\ell_j \geq 0 
\] 
for all $d \in \N$, 
then we define the function $f_{j+1}$ as follows: 
\[
f_{j+1}(k) = 
\begin{cases}
f_{j}(k) - q_j\ell_j & \text{if $k  \equiv 0 \pmod{\ell_j}$} \\
f_j(k) & \text{if $k \not\equiv 0 \pmod{\ell_j}$.} 
\end{cases}
\]
Note that $f_{j+1}(k) = 0$ for all $k \leq \ell_j$. 
\eenum
\eenum
\eenum

If $f_j$ is the constant function -1 for some $j$, then $f$ is not the fixed point 
counting function of a permutation.  

If $f_j$ is the zero function for some $j$, and if $r$ is the smallest positive integer 
such that $f_r$  is the zero function, 
then there are finite sequences of positive 
integers $(\ell_j)_{j=1}^r$ and $(q_j)_{j=1}^r$ such that 
$(\ell_j)_{j=1}^r$  is strictly increasing and $f$ is the fixed point 
counting function of a permutation $\sigma$ of a nonempty finite set $X$ with   
\[
|X| =\sum_{j=1}^r q_j\ell_j \qqand C_{\sigma}(\ell_j) = q_j \quad\text{for $j \in \{1,\ldots, r\}$.}
\] 

If $f_j$ is neither the zero function nor the constant function -1 for all $j \in \N$, then 
$(f_j)_{j=1}^{\infty}$ is an infinite sequence of nonzero nonnegative 
arithmetic functions. 
There are infinite sequences of positive 
integers $(\ell_j)_{j=1}^{\infty}$ and $(q_j)_{j=1}^{\infty}$ such that 
$(\ell_j)_{j=1}^{\infty}$  is strictly increasing and $f$ is the fixed point 
counting function of a permutation $\sigma$ of a nonempty finite set $X$ with   
\[
C_{\sigma}(\ell_j) = q_j\ell_j \qquad\text{for all $j \in \N$} 
\]
and 
\[
C_{\sigma}(k) = 0\qquad\text{for all $k \in \N$ such that $k \neq \ell_j$ for all $j \in \N$.}
\]
Thus, the algorithm determines if $f$ is the fixed point counting function of a permutation 
and, if so, constructs a set $X$ and a permutation $\sigma \in \Perm(X)$ such that 
$f(k) = F_X(\sigma^k)$ for all $k \in \N$.

\section{Periodic fixed point counting functions}

We shall determine the permutations with periodic fixed point counting functions.
Here is a simple example. 
If $\sigma$ has order $m$ in the group $\Perm(X)$, 
then for all $k \in \Z$ we have $\sigma^{k+m} = \sigma^k$.  Thus,  
\beq                \label{FixedPoint:order=period}
F_X\left( \sigma^{k+m} \right) = F_X\left( \sigma^k \right)
\eeq
and the fixed point counting function $F_X(\sigma^k)$ is periodic with period $m$. 
In particular, the fixed point counting function of a permutation of a finite set is periodic.

\bt                                             \label{FixedPoint:theorem:periodic} 
 Let $X$ be a nonempty set and let $\sigma \in \Perm(X)$.  
 The fixed point counting  function $F_X(\sigma^k)$ 
is periodic if and only if the permutation $\sigma$ is bounded.   
\et

\begin{proof}
Let $\sigma \in \Perm(X)$ and let $\hat{X}_0$ and $\hat{X}_{\infty}$  
be the sets defined by~\eqref{FixedPoint:X-0} and~\eqref{FixedPoint:X-infty}. 
If $\hat{X}_0 = \emptyset$, then $X = \hat{X}_{\infty}$ 
and  $F_X(\sigma^k) $ is the periodic constant function 0.
The permutation $\sigma$  is of bounded type because it contains no cycle of finite length.  

If $\hat{X}_0 \neq \emptyset$, then, by Theorem~\ref{FixedPoint:theorem:PartitionFinite}, 
$F_X(\sigma^k)  = F_{\hat{X}_0} (\hat{\sigma_0}^k)$, 
where $\hat{\sigma}_0$ is the restriction of $\sigma$ to $\hat{X}_0$.
Thus, we can assume without loss of generality 
that $X = \hat{X}_0$ and the permutation $\sigma$ is of finite type. 

Let $\sigma = \prod_{i \in I} \gamma_i$ be the cyclic decomposition of $\sigma$, where 
$\gamma_i$ is a cycle of finite length $\ell_i$ on the set $X_i = \{x_{i,1},\ldots, x_{i,\ell_i} \}$.  
Let $L = \{\ell_i: i \in I\}$ be the set of lengths of the cycles in $\sigma$.  
The permutation $\sigma$ is bounded if and only if  $L$ is a finite set of positive integers.  
Suppose that $L$ is finite and that   $m$ is the least common multiple of the integers in $L$.  
For all $i \in I$ and $k \in \N$ we have $\gamma_i^m = \id_X$ and $\gamma_i^{k+m} = \gamma_i^k$.
This implies $\sigma^{k+m} = \sigma^k$ and so  
$F_X\left( \sigma^{k+m} \right) = F_X(\sigma^k)$ for all $k \in \Z$. 
Thus,  the fixed point counting function is periodic with period $m$.  

If the set $L = \{\ell_i: i \in I\}$ of lengths of   cycles in $\sigma$ is unbounded, 
then there is a strictly increasing infinite sequence $\Lambda = (\ell_{i_j})_{j=1}^{\infty}$ 
of elements of $L$.  Let $m_r = \lcm(\ell_{i_1},\ldots, \ell_{i_r} )$ 
be the least common multiple of the first $r$ integers in the sequence $\Lambda$. 
For all $j \in \{1,2,\ldots, r\}$ we have  
\[
\lambda_{i_j}^{m_r} = \lambda_{i_j}^{\ell_{i_j}} = \id_X
\]
and so 
\[
\sigma^{m_r}(x) = x
\]
for all 
\[
x \in \bigcup_{j=1}^r X_{i_j}.  
\]
It follows that 
\[
F_X\left( \sigma^{m_r} \right) \geq  \sum_{j=1}^r |X_{i_j}| = \sum_{j=1}^r \ell_{i_j}
\]
and 
\[
\lim_{ r \rightarrow \infty} F_X\left( \sigma^{m_r} \right) = \infty. 
\]
The function $F_X(\sigma^k)$ is unbounded and, therefore, not periodic. 
This completes the proof. 
\end{proof}

\end{document}